\documentstyle[12pt,amscd]{amsart}
\newtheorem{thm}[equation]{Theorem}
\numberwithin{equation}{section}
\newtheorem{cor}[equation]{Corollary}
\newtheorem{expl}[equation]{Example}
\newtheorem{rmk}[equation]{Remark}

\newtheorem{lem}[equation]{Lemma}

\newtheorem{defin}[equation]{Definition}

\newtheorem{prop}[equation]{Proposition}

\begin{document}
\raggedbottom \voffset=-.7truein \hoffset=0truein \vsize=8truein
\hsize=6truein \textheight=8truein \textwidth=6truein
\baselineskip=18truept
\def\mapright#1{\ \smash{\mathop{\longrightarrow}\limits^{#1}}\ }
\def\mapleft#1{\smash{\mathop{\longleftarrow}\limits^{#1}}}
\def\mapup#1{\Big\uparrow\rlap{$\vcenter {\hbox {$#1$}}$}}
\def\mapdown#1{\Big\downarrow\rlap{$\vcenter {\hbox {$\ssize{#1}$}}$}}
\def\mapne#1{\nearrow\rlap{$\vcenter {\hbox {$#1$}}$}}
\def\mapse#1{\searrow\rlap{$\vcenter {\hbox {$\ssize{#1}$}}$}}
\def\mapr#1{\smash{\mathop{\rightarrow}\limits^{#1}}}
\def\l{\lambda}
\def\TK{T_{K/2}}
\def\ss{\smallskip}
\def\vp{v_1^{-1}\pi}
\def\at{{\widetilde\alpha}}
\def\sm{\wedge}
\def\la{\langle}
\def\ra{\rangle}
\def\on{\operatorname}
\def\spin{\on{Spin}}
\def\kbar{{\overline k}}
\def\lbar{{\overline{\ell}}}
\def\qed{\quad\rule{8pt}{8pt}\bigskip}
\def\ssize{\scriptstyle}
\def\a{\alpha}
\def\Z{{\Bbb Z}}
\def\im{\on{im}}
\def\ct{\widetilde{C}}
\def\gdt{\widetilde{\on{gd}}}
\def\ext{\on{Ext}}
\def\sq{\on{Sq}}
\def\eps{\epsilon}
\def\ar#1{\stackrel {#1}{\rightarrow}}
\def\br{{\bold R}}
\def\bc{{\bold C}}
\def\bh{{\bold H}}
\def\bx{{\bold x}}
\def\si{\sigma}
\def\Ebar{{\overline E}}
\def\dbar{{\overline d}}
\def\Sum{\sum}
\def\tfrac{\textstyle\frac}
\def\tb{\textstyle\binom}
\def\Si{\Sigma}
\def\w{\wedge}
\def\equ{\begin{equation}}
\def\b{\beta}
\def\G{\Gamma}
\def\g{\gamma}
\def\psit{\widetilde{\Psi}}
\def\tht{\widetilde{\Theta}}
\def\psiu{{\underline{\Psi}}}
\def\thu{{\underline{\Theta}}}
\def\aee{A_{\text{ee}}}
\def\aeo{A_{\text{eo}}}
\def\aoo{A_{\text{oo}}}
\def\aoe{A_{\text{oe}}}
\def\fbar{{\overline f}}
\def\endeq{\end{equation}}
\def\sn{S^{2n+1}}
\def\zp{\bold Z_p}
\def\A{{\cal A}}
\def\P{{\cal P}}
\def\cj{{\cal J}}
\def\zt{{\bold Z}_2}
\def\bs{{\bold s}}
\def\bof{{\bold f}}
\def\bq{{\bold Q}}
\def\be{{\bold e}}
\def\Hom{\on{Hom}}
\def\ker{\on{ker}}
\def\coker{\on{coker}}
\def\da{\downarrow}
\def\colim{\operatornamewithlimits{colim}}
\def\zphat{\bz_2^\wedge}
\def\io{\iota}
\def\Om{\Omega}
\def\u{{\cal U}}
\def\e{{\cal E}}
\def\exp{\on{exp}}
\def\wbar{{\overline w}}
\def\xbar{{\overline x}}
\def\ybar{{\overline y}}
\def\zbar{{\overline z}}
\def\ebar{{\overline e}}
\def\nbar{{\overline n}}
\def\rbar{{\overline r}}
\def\et{{\widetilde E}}
\def\ni{\noindent}
\def\coef{\on{coef}}
\def\den{\on{den}}
\def\lcm{\on{l.c.m.}}
\def\vi{v_1^{-1}}
\def\ot{\otimes}
\def\psibar{{\overline\psi}}
\def\mhat{{\hat m}}
\def\exc{\on{exc}}
\def\ms{\medskip}
\def\ehat{{\hat e}}
\def\etao{{\eta_{\text{od}}}}
\def\etae{{\eta_{\text{ev}}}}
\def\dirlim{\operatornamewithlimits{dirlim}}
\def\gt{\widetilde{L}}
\def\lt{\widetilde{\lambda}}
\def\st{\widetilde{s}}
\def\ft{\widetilde{f}}
\def\sgd{\on{sgd}}
\def\lfl{\lfloor}
\def\rfl{\rfloor}
\def\ord{\on{ord}}
\def\gd{{\on{gd}}}
\def\rk{{{\on{rk}}_2}}
\def\nbar{{\overline{n}}}
\def\lg{{\on{lg}}}
\def\N{{\Bbb N}}
\def\Z{{\Bbb Z}}
\def\Q{{\Bbb Q}}
\def\R{{\Bbb R}}
\def\C{{\Bbb C}}
\def\mo{\on{mod}}
\def\vexp{v_1^{-1}\exp}
\def\Remark{\noindent{\it  Remark}}
\title[$v_1$-periodic homotopy groups of $DI(4)$]
{$v_1$-periodic homotopy groups of the Dwyer-Wilkerson space}
\author{Martin Bendersky}
\address{Hunter College, CUNY\\NY, NY 01220}
\email{mbenders@@hunter.cuny.edu}
\author{Donald M. Davis}
\address{Lehigh University\\Bethlehem, PA 18015}
\email{dmd1@@lehigh.edu}

\date{June 7, 2007}
\keywords{$v_1$-periodic homotopy groups, $p$-compact groups, Adams operations, $K$-theory}
\subjclass[2000] {55Q52, 57T20, 55N15.}

\maketitle
\begin{abstract} The Dwyer-Wilkerson space $DI(4)$ is the only exotic
2-compact group. We compute its $v_1$-periodic homotopy groups $\vp_*(DI(4))$.
\end{abstract}

\section{Introduction}\label{intro}
In \cite{DW1}, Dwyer and Wilkerson constructed a $2$-complete space $BDI(4)$,
so named because its $F_2$-cohomology groups form an algebra isomorphic to the
ring of rank-4 mod-2 Dickson invariants. Its loop space, called $DI(4)$, has
$H^*(DI(4);F_2)$ finite. In \cite{DW2}, they then defined a $p$-compact group to be
a pair $(X,BX)$, such that $X=\Omega BX$ (hence $X$ is redundant), $BX$ is connected
and $p$-complete, and $H^*(X;F_p)$ is finite. In \cite{Gro}, Andersen and Grodal proved that
$(DI(4), BDI(4))$ is the only simple 2-compact group not arising as the 2-completion of
a compact connected Lie group.

The $p$-primary $v_1$-periodic homotopy groups of a topological space $X$, defined in \cite{DM}
and denoted
$\vp_*(X)_{(p)}$ or just $\vp_*(X)$ if the prime is clear, are a first approximation
to the $p$-primary homotopy groups. Roughly, they are a localization of the portion of the
actual homotopy groups detected by $p$-local $K$-theory. In \cite{DHHA}, the second author
completed a 13-year project, often in collaboration with the first author, of determining
$\vp_*(X)_{(p)}$ for all compact simple Lie groups and all primes $p$

In this paper, we determine the 2-primary groups $\vp_*(DI(4))$. Here and throughout,
$\nu(-)$ denotes the exponent of 2 in the prime factorization of an integer.
\begin{thm}\label{main} 
For any integer $i$, let $e_i=\min(21,4+\nu(i-90627))$. Then
$$\vp_{8i+d}(DI(4))\approx\begin{cases}\Z/2^{e_i}\oplus\Z/2&d=1\\
\Z/2^{e_i}&d=2\\
0&d=3,4\\
\Z/8&d=5\\
\Z/8\oplus\Z/2&d=6\\
\Z/2\oplus \Z/2\oplus\Z/2&d=7,8.\end{cases}$$
\end{thm}

Since every $v_1$-periodic homotopy group is a subgroup of some actual
homotopy group, this result implies that $\exp_2(DI(4))\ge21$, i.e., some homotopy group
of $DI(4)$ has an element of order $2^{21}$. It would be interesting to know
whether this bound is sharp.

Our proof involves studying the spectrum $\Phi_1DI(4)$
which satisfies $\pi_*(\Phi_1DI(4))\approx\vp_*(DI(4))$. We will relate $\Phi_1DI(4))$
to the 2-completed $K$-theoretic pseudosphere $T_{K/2}$ discussed in \cite[8.6]{BoTop}.
We will prove the following surprising result, which was pointed out by Pete Bousfield.
\begin{thm}\label{spherethm} There is an equivalence of spectra $$\Phi_1DI(4)\simeq\Sigma^{725019}T_{K/2}\w M(2^{21}),$$
where $M(2^{21})$ is a mod $2^{21}$ Moore spectrum.
\end{thm} 
In Section \ref{psspsec}, we will give the easy deduction of
Theorem \ref{main} from Theorem \ref{spherethm}.
As an immediate corollary of \ref{spherethm}, we deduce that the $2^{21}$ bound on
$\pi_*(\Phi_1DI(4))$ is induced from a bound on the spectrum itself.
\begin{cor} The exponent of the spectrum $\Phi_1DI(4)$ is $2^{21}$; i.e.,
$2^{e}1_{\Phi_1DI(4)}$ is null if and only if $e\ge21$.\end{cor}

In \cite{Bo}, Bousfield presented a framework that enables determination of the
$v_1$-periodic homotopy groups of many simply-connected $H$-spaces $X$ from their
united $K$-theory groups and Adams operations. The intermediate step is
$KO^*(\Phi_1 X)$. (All of our $K^*(-)$ and $KO^*(-)$-groups have coefficients in the 2-adic integers $\hat \Z_2$, which we omit
from our notation.) Our first proof of Theorem \ref{main} used Bousfield's exact sequence
\cite[9.2]{Bo} which relates $\vp_*(X)$ with $KO^*(\Phi_1X)$, but the approach via the
pseudosphere, which we present here, is stronger and more elegant. The insight for Theorem \ref{spherethm} is the observation that
the two spectra have isomorphic Adams modules $KO^*(-)$.

In several earlier e-mails,
Bousfield explained to the authors how the results of \cite{Bo} should
enable us to determine $KO^*(\Phi_1DI(4))$. In Section \ref{Bosec}, we present our account
of these ideas of Bousfield.  We thank him profusely for sharing his insights
with us.

The other main input is the Adams operations in $K^*(BDI(4))$.  
In \cite{OS}, Osse and Suter showed that $K^*(BDI(4))$ is a power series algebra
on three specific generators, and gave some information toward the determination of the
Adams operations. In private communication in 2005, Suter expanded on this to give
explicit formulas for $\psi^k$ in $K^*(BDI(4))$. We are very grateful to him for
sharing this information. In Section \ref{compsec}, we will
explain these calculations and also how they then lead to the determination
of $KO^*(\Phi_1DI(4))$.

\section{Adams operations}\label{compsec}
In this section, we present Suter's determination of $\psi^k$ in $K^*(BDI(4))$ and
state a result, proved in Section \ref{Bosec}, that allows us to determine
$KO^*(\Phi_1DI(4))$ from these Adams operations.

Our first result, communicated by Suter, is the following determination of Adams operations
in $K^*(BDI(4))$. An element of $K^*(X)$ is called {\it real} if it is in the image of
$KO^*(X)@>c>> K^*(X)$.
\begin{thm}\label{psik} (Suter) There is an isomorphism of algebras
\begin{equation}\label{three}K^*(BDI(4))\approx \hat\Z_2[\![\xi_8,\xi_{12},\xi_{24}]\!]\end{equation}
such that the generators are in $K^0(-)$ and are real, $\psi^{-1}=1$, and the matrices of $\psi^2$ and $\psi^3$
on the three generators, mod decomposables, are
$$\Psi^2\equiv\begin{pmatrix}2^4&0&0\\ -2&2^6&0\\ 0&-2&2^{14}\end{pmatrix},\qquad\Psi^3\equiv
\begin{pmatrix}3^4&0&0\\-3^3&3^6&0\\ 36/527&-3^5\cdot41/17&3^{14}\end{pmatrix}.$$
\end{thm}
\begin{pf} The subscripts of the generators indicates their ``filtration," meaning
the dimension of the smallest skeleton on which they are nontrivial. A standard
property of Adams operations is that if $\xi$ has filtration $2r$, then $\psi^k(\xi)$
equals $k^r\xi$ plus elements of higher filtration.

The isomorphism (\ref{three}) is derived in \cite[p.184]{OS} along with the additional information
that $4\xi_{24}-\xi_{12}^2$ has filtration 28, and \begin{eqnarray}\label{l2}\xi_{12}&=&\lambda^2(\xi_8)+8\xi_8\\ \nonumber\xi_{24}&=&
\lambda^2(\xi_{12})+32\xi_{12}+c_1\xi_8^2+c_2\xi_8^3+c_3\xi_8\xi_{12},\end{eqnarray}
for certain explicit even coefficients $c_i$. 

The Atiyah-Hirzebruch spectral sequence easily shows that $\xi_8$ is real, since the 11-skeleton
of $BDI(4)$ equals $S^8$. Since $\lambda^2(c(\theta))=c(\lambda^2(\theta))$, and products of real
bundles are real, we deduce from (\ref{l2}) that $\xi_{12}$ and $\xi_{24}$ are also real. Since $tc=c$, where $t$
denotes conjugation, which corresponds to $\psi^{-1}$, we obtain that the generators are invariant under
$\psi^{-1}$, and hence so is all of $K^*(BDI(4))$.

We compute Adams operations mod decomposables, writing $\equiv$ for equivalence mod decomposables.
Because $4\xi_{24}-\xi_{12}^2$ has filtration 28, we obtain \begin{equation}\label{strange}\psi^k(\xi_{24})\equiv k^{14}\xi_{24}.\end{equation}
Here we use, from \cite[p.183]{OS}, that all elements of $K^*(BDI(4))$ of filtration greater than 28
are decomposable. Equation (\ref{strange}) may seem surprising, since $\xi_{24}$ has filtration 24,
but there is a class $\xi_{28}$ such that $4\xi_{24}-\xi_{12}^2=\xi_{28}$, and we can have
$\psi^k(\xi_{24})\equiv k^{12}\xi_{24}+\a_k\xi_{28}$ consistently with (\ref{strange}).

Using (\ref{l2}) and that $\psi^2\equiv -2\lambda^2$ mod decomposables, we obtain
\begin{eqnarray}\label{psi2}\psi^2(\xi_8)&\equiv&2^4\xi_8-2\xi_{12}\\
\psi^2(\xi_{12})&\equiv&2^6\xi_{12}-2\xi_{24}\nonumber,\end{eqnarray}
yielding the matrix $\Psi^2$ in the theorem.

Applying $\psi^2\psi^3=\psi^3\psi^2$ to $\psi^3(\xi_{12})\equiv3^6\xi_{12}+\gamma\xi_{24}$
yields $-2\cdot3^6+2^{14}\gamma=2^6\gamma-2\cdot3^{14}$, from which we obtain $\gamma=-3^5\cdot41/17$.
Applying the same relation to $\psi^3(\xi_8)=3^4\xi_8+\a\xi_{12}+\b\xi_{24}$, coefficients
of $\xi_{12}$ yield $-2\cdot3^4+\a\cdot2^6=2^4\a-2\cdot3^6$ and hence $\a=-3^3$.
Now coefficients of $\xi_{24}$ yield $-2\a+2^{14}\b=2^4\b-2\gamma$ and hence $\b=36/527$.
\end{pf}

Let $\Phi_1(-)$ denote the functor from spaces to $K/2_*$-local spectra described in
\cite[9.1]{Bo}, which satisfies $\vp_*X\approx\pi_*\tau_2\Phi_1X$, where $\tau_2\Phi_1X$ is the 2-torsion 
part of $\Phi_1X$. 
 In Section \ref{Bosec}, we will use results of Bousfield in \cite{Bo} to prove the following result. 
Aspects of Theorem \ref{psik}, such as $K^*(BDI(4))$ being a power
series algebra on real generators, are also used in proving this theorem.

Recall that $KO^*(-)$ has period 8.
\begin{thm}\label{KO} The groups $KO^i(\Phi_1DI(4))$ are $0$ if $i\equiv0,1,2\mod 8$, and $K^0(\Phi_1 DI(4))=0$. Let $M$ denote a
free $\hat\Z_2$-module on three generators, acted on by $\psi^2$ and $\psi^3$ by the matrices of Theorem
\ref{psik}, with $\psi^{-1}=1$. Let $\theta=\frac12\psi^2$ act on $M$. Then there are exact sequences
\begin{eqnarray*}0&\to& 2M@>\theta>>2M\to KO^3(\Phi_1DI(4))\to0\to0\to KO^4(\Phi_1DI(4))\to M/2\\
&@>\theta>>& M/2\to KO^5(\Phi_1DI(4))\to M/2@>\theta>>M/2\to KO^6(\Phi_1DI(4))\to M\\
&@>\theta>>& M\to KO^7(\Phi_1DI(4))\to 0
\end{eqnarray*}
and
$$0\to M@>\theta>> M\to K^1(\Phi DI(4))\to 0.$$
For $k=-1$ and $3$, the action of $\psi^k$ in $KO^{2j-1}(\Phi_1DI(4))$, $KO^{2j-2}(\Phi_1DI(4))$, and $K^{2j-1}(\Phi DI(4))$ agrees with $k^{-j}\psi^k$ in adjacent $M$-terms.
\end{thm}

In the remainder of this section, we use \ref{psik} and \ref{KO} to give explicit formulas for the Adams module $KO^i(\Phi_1 DI(4))$. A similar argument works for $K^*(\Phi_1 DI(4))$.
If $g_1$, $g_2$, and $g_3$ denote the three generators of $M$, then the action of $\theta$ is given by
\begin{eqnarray*}\theta(g_1)&=&8g_1-g_2\\
\theta(g_2)&=&2^5g_2-g_3\\
\theta(g_3)&=&2^{13}g_3.\end{eqnarray*}
Clearly $\theta$ is injective on $M$ and $2M$. We have $KO^7(\Phi_1DI(4))\approx\coker(\theta|M)\approx\Z/2^{21}$ with generator $g_1$; note that $g_2=2^3g_1$ in this cokernel, and then $g_3=2^8g_1$. Similarly $KO^3(\Phi_1DI(4))\approx\coker(\theta|2M)\approx\Z/2^{21}$.
Also $KO^4(\Phi_1DI(4))\approx\ker(\theta|M/2)=\Z/2$ with generator $g_3$, while
$KO^6(\Phi_1DI(4))\approx\coker(\theta|M/2)=\Z/2$ with generator $g_1$. There is a short exact sequence
$$0\to\coker(\theta|M/2)\to KO^5(\Phi_1DI(4))\to\ker(\theta|M/2)\to 0,$$
with the groups at either end being $\Z/2$ as before. To see that this short exact sequence is split,
we use the map $S^7@>f>> DI(4)$ which is inclusion of the bottom cell. The morphism $f^*$ sends the first
summand of $KO^5(\Phi_1DI(4))$ to one of the two $\Z/2$-summands of $KO^5(\Phi_1S^7)$, providing a splitting
homomorphism.
 Thus we have proved the first part of the following result.
\begin{thm}\label{gps} We have
\begin{eqnarray*}K^i(\Phi_1DI(4))&\approx&\begin{cases}0&i=0\\ \Z/2^{21}&i=1,\end{cases}\\
KO^i(\Phi_1DI(4))&\approx&\begin{cases}0&i=0,1,2\\
\Z/2^{21}&i=3,7\\
\Z/2&i=4,6\\
\Z/2\oplus\Z/2&i=5.\end{cases}\end{eqnarray*}
For $k=-1$ and $3$, we have $\psi^k=1$ on the $\Z/2$'s, and on $KO^{2j-1}(\Phi_1DI(4))$ with $j$ even and $K^{2j-1}(\Phi_1DI(4))$, $\psi^{-1}=(-1)^j$ and
$$\psi^3=3^{-j}(3^4-3^3\cdot2^3+\tfrac{36}{527}2^8).$$\end{thm}

\begin{pf*}{Completion of proof}
To obtain $\psi^3$ on the $\Z/2$'s, we use the last part of Theorem \ref{KO} and the matrix $\Psi^3$ of
Theorem \ref{psik}. If $\psi^3$ is as in $\Psi^3$, then, mod 2, $\psi^3-1$ sends
$g_1\mapsto g_2$, $g_2\mapsto g_3$, and $g_3\mapsto0$. Thus  $\psi^3-1$ equals 0 on $KO^4(\Phi_1DI(4))$ and $KO^6(\Phi_1DI(4))$. Clearly $\psi^{-1}=1$ on these groups.

To see that $\psi^k-1$ is 0 on $KO^5(\Phi_1DI(4))$, we use the commutative diagram
$$\begin{CD} 0@>>>\Z/2@>i>> KO^5(\Phi_1DI(4)) @>\rho>>\Z/2@>>>0\\
@. @V\approx VV @V f^* VV @V0VV @.\\
0@>>>\Z/2 @>>> KO^5(\Phi_1S^7)@>>>\Z/2@>>>0.\end{CD}$$
We can choose generators $G_1$ and $G_2$ of $KO^5(\Phi_1DI(4))\approx\Z/2\oplus\Z/2$ so that $G_1\in\im(i)$,
$\rho(G_2)\ne0$, and $f^*(G_2)=0$. Since $\psi^k-1=0$ on the $\Z/2$'s on either side of $KO^5(\Phi_1DI(4))$,
the only way that $\psi^k-1$ could be nonzero on $KO^5(\Phi_1DI(4))$ is if $(\psi^k-1)(G_2)=G_1$. However this yields the
contradiction
$$0=(\psi^k-1)f^*G_2=f^*(\psi^k-1)G_2=f^*(G_1)\ne0.$$

On $KO^{2j-1}(\Phi_1DI(4))$ with $j$ even and $K^{2j-1}(\Phi_1DI(4))$, $\psi^3$ sends the generator $g_1$ to
$$3^{-j}(3^4g_1-3^3g_2+\tfrac{36}{527}g_3)=3^{-j}(3^4-3^3\cdot2^3+\tfrac{36}{527}2^8)g_1,$$ 
and $\psi^{-1}(g_1)=(-1)^jg_1$ by Theorem \ref{KO}.
\end{pf*}

\section{Relationship with pseudosphere}\label{psspsec}
In this section, we prove Theorems \ref{spherethm} and \ref{main}.

Following \cite[8.6]{BoTop}, we let $T=S^0\cup_\eta e^2\cup_2 e^3$, and consider its
$K/2$-localization $T_{K/2}$. The groups $\pi_*(T_{K/2})$ are given in \cite[8.8]{BoTop}, while 
the Adams module  is given by
\begin{eqnarray*}K^i(\TK)&=&\begin{cases}\hat\Z_2&i\text{ even, with }\psi^k=k^{-i/2}\\
0&i\text{ odd;}\end{cases}\\
KO^i(\TK)&=&\begin{cases}\hat\Z_2&i\equiv0\mod4,\text{ with }\psi^k=k^{-i/2}\\
\Z/2&i=2,3,\text{ with }\psi^k=1\\
0&i=1,5,6,7.\end{cases}
\end{eqnarray*}
Bousfield calls this the 2-completed $K$-theoretic pseudosphere. Closely related spectra
have been also considered in \cite{HMS} and \cite{BDM}.

Let $M(n)=S^{-1}\cup_n e^0$ denote the mod $n$ Moore spectrum. Then, for $e>1$ and $k$ odd,
\begin{eqnarray}\nonumber K^i(\TK\w M(2^e))&=&\begin{cases}\Z/2^e&i\text{ even, with }\psi^k=k^{-i/2}\\
0&i\text{ odd;}\end{cases}\\
\label{KY} KO^i(\TK\w M(2^e))&=&\begin{cases}\Z/2^e&i\equiv0\mod4,\text{ with }\psi^k=k^{-i/2}\\
\Z/2&i=1,3,\text{ with }\psi^k=1\\
\Z/2\oplus\Z/2&i=2,\text{ with }\psi^k=1\\
0&i=5,6,7.\end{cases}
\end{eqnarray}
\begin{pf} Let $Y=\TK\w M(2^e)$. Most of (\ref{KY}) is immediate from the exact sequence
$$@>2^e>>KO^i(\TK)@>>> KO^i(Y)@>>> KO^{i+1}(\TK)@>2^e>>.$$
To see that $KO^2(Y)=\Z/2\oplus\Z/2$ and not $\Z/4$, one can first note that
\begin{equation}\label{Msplit}M(2^e)\w M(2)\simeq \Sigma^{-1}M(2)\vee M(2).\end{equation}
The exact sequence
\begin{equation}\label{KOY}KO^2(Y)@>2>>KO^2(Y)@>>>KO^2(Y\w M(2))@>>>KO^3(Y)@>2>>\end{equation}
implies that if $KO^2(Y)=\Z/4$, then $|KO^2(Y\w M(2))|=4$. However, by (\ref{Msplit}),
\begin{equation}\label{KM}KO^2(Y\w M(2))\approx KO^2(\TK\w M(2))\oplus KO^3(\TK\w M(2)).\end{equation}
Also, there is a cofiber sequence
\begin{equation}\label{A1cof}T\w M(2)\to \Sigma^{-1}A_1 \to \Sigma^5 M(2),\end{equation}
where $H^*(A_1;F_2)$ is isomorphic to the subalgebra of the mod 2 Steenrod algebra generated
by $\sq^1$ and $\sq^2$, and satisfies $KO^*(A_1)=0$. Thus
$$KO^i(\TK\w M(2))\approx KO^i(\Sigma^4M(2))\approx\begin{cases}\Z/4&i=2\\
\Z/2&i=3,\end{cases}$$
so that $|KO^2(Y\w M(2))|=8$, contradicting a consequence of the hypothesis that $KO^2(Y)=\Z/4$.

We conclude the proof by showing that, for odd $k$, $\psi^k=1$ on $KO^2(Y)$. First note that $\psi^k=1$
on $KO^*(M(2))$. This follows immediately from the Adams operations on the sphere, except for $\psi^k$ on $KO^{-2}(M(2))\approx\Z/4$. This is isomorphic
to $\widetilde{KO}(RP^2)$, where $\psi^k=1$ is well-known. Now use (\ref{MTM}) to deduce that $\psi^k=1$ on
$KO^*(T_{K/2}\w M(2))$, and then (\ref{KM}) to see that $\psi^k=1$ on $KO^2(Y\w M(2))$. Finally, use (\ref{KOY})
to deduce that $\psi^k=1$ on $KO^2(Y)$.
\end{pf}

Comparison of \ref{gps} and (\ref{KY}) yields an isomorphism of graded abelian groups 
\begin{equation}\label{KOiso}KO^*(\Sigma^{8L+3}T_{K/2}\w M(2^{21}))\approx KO^*(\Phi_1DI(4))\end{equation}
for any integer $L$. We will show that if $L=90627$, then the Adams operations
agree too. By \cite[6.4]{BoKs}, it suffices to prove they agree for $\psi^3$ and $\psi^{-1}$.

Note that one way of distinguishing a $K$-theoretic pseudosphere from
a sphere is that in $KO^*(\text{sphere})$ (resp. $KO^*(\text{pseudosphere})$)
the $\Z/2$-groups are in dimensions 1 and 2 less than the dimensions in which
$\psi^3\equiv1$ mod 16 (resp. $\psi^3\equiv9$ mod 16), and similarly after smashing
with a mod $2^e$ Moore spectrum. Since
$3^4-6^3+\frac{36}{527}2^8\equiv9$ mod 16, the $\Z/2$-groups in $KO^*(\Phi_1DI(4))$
occur in dimensions 1, 2, and 3 less that the dimension in which $\psi^3\equiv9$ mod 16, and so 
$\Phi_1DI(4)$ should be identified with a suspension of $T_{K/2}\w M(2^{21})$ and not
$S_{K/2}\w M(2^{21})$.

In $KO^{4t-1}(\Sigma^{8L+3}T_{K/2}\w M(2^{21}))$, $\psi^3=3^{-2(t-2L-1)}$ and $\psi^{-1}=1$.
Thus if $L$ satisfies
\begin{equation}\label{Leq} 3^{4L+2}\equiv 3^4-6^3+\tfrac{36}{527}2^8\mod 2^{21},\end{equation}
then $KO^*(\Sigma^{8L+3}T_{K/2}\w M(2^{21}))$ and $KO^*(\Phi_1 DI(4))$ will be isomorphic
Adams modules. {\tt Maple} easily verifies that (\ref{Leq}) is satisfied for $L=90627$.

A way in which this number $L$ can be found begins with the mod $2^{18}$ equation
$$\sum_{i=1}^6\tbinom{2L-1}i 8^{i-1}\equiv\frac{3^{4L-2}-1}8\equiv\tfrac19(\tfrac{2^7}{527}-3)\equiv192725,$$
where we use {\tt Maple} at the last step. This easily implies $L\equiv3$ mod 8, and so we let $L=8b+3$.
Again using {\tt Maple} and working mod $2^{18}$ we compute
$$\sum_{i=1}^6\tbinom{16b+5}i8^{i-1}-192725\equiv 2^{10}u_0+2^4u_1b+2^{10}u_2b^2+2^{17}b^3,$$
with $u_i$ odd. Thus we must have $b\equiv64$ mod 128, and so $L\equiv515$ mod $2^{10}$.
Several more steps of this type lead to the desired value of $L$.

Thus, in the terminology of \ref{CRexpl}, we have proved the following result.
\begin{prop} If $L=90627$, then there is an isomorphism of Adams modules
$$K_{CR}^*(\Sigma^{8L+3}T_{K/2}\w M(2^{21}))\approx K^*_{CR}(\Phi_1 DI(4)).$$
\end{prop}
Theorem \ref{spherethm} follows immediately from this using the remarkable
\cite[5.3]{BoTop}, which says, among other things, that 2-local spectra $X$
having some $K^i(X)=0$ are determined up to equivalence by their Adams module $K_{CR}^*(X)$.
Theorem \ref{main} follows immediately from Theorem \ref{spherethm} and the 
following result.
\begin{prop}\label{piT} For all integers $i$,
$$\pi_{8i+d}(T_{K/2}\w M(2^{21}))\approx\begin{cases}\Z/2\oplus\Z/2^{\min(21,\nu(i)+4)}&d=-2\\
\Z/2^{\min(21,\nu(i)+4)}&d=-1\\
0&d=0,1\\
\Z/8&d=2\\\Z/2\oplus\Z/8&d=3\\
\Z/2\oplus\Z/2\oplus\Z/2&d=4,5.\end{cases}$$
\end{prop}
\begin{pf} For the most part, these groups are immediate from the groups $\pi_*(T_{K/2})$ given in 
\cite[8.8]{BoTop} and the exact sequence
\begin{equation}\label{les}@>2^{21}>>\pi_{j+1}(T_{K/2})\to \pi_j(T_{K/2}\w M(2^{21}))\to \pi_j(T_{K/2})@>2^{21}>>.
\end{equation}
All that needs to be done is to show that the following short exact sequences, obtained from (\ref{les}), are split.
\begin{eqnarray}\label{4seq}0\ \to\ \Z/2\ \to&\pi_{8i+3}(T_{K/2}\w M(2^{21}))&\to\ \Z/8\ \to\ 0\\
0\ \to\ \Z/2\oplus\Z/2\ \to&\pi_{8i+4}(T_{K/2}\w M(2^{21}))&\to\ \Z/2\ \to\ 0\nonumber\\
0\ \to\ \Z/2\ \to&\pi_{8i+5}(T_{K/2}\w M(2^{21}))&\to\ \Z/2\oplus\Z/2\ \to\ 0\nonumber\\
0\ \to\ \Z/2^{\min(21,\nu(i)+4)}\ \to&\pi_{8i-2}(T_{K/2}\w M(2^{21}))&\to\ \Z/2\ \to\ 0.\nonumber\end{eqnarray}
Let $Y=T_{K/2}\w M(2^{21})$. We consider the exact sequence for $\pi_*(Y\w M(2))$,
\begin{equation}\label{YMseq}@>2>>\pi_{i+1}(Y)\to \pi_i(Y\w M(2))\to \pi_i(Y)@>2>>.\end{equation}
If the four sequences (\ref{4seq}) are all split, then by (\ref{YMseq}) the groups $\pi_{8i+d}(Y\w M(2))$
for $d=2,3,4,5,-2$ have orders $2^3$, $2^5$, $2^6$, $2^5$, and $2^3$, respectively, but if any of the sequences
(\ref{4seq}) fails to split, then some of the orders $|\pi_{8i+d}(Y\w M(2))|$ will have values smaller than those
listed here.

By (\ref{Msplit}),
$$\pi_i(Y\w M(2))\approx\pi_{i+1}(T_{K/2}\w M(2))\oplus \pi_i(T_{K/2}\w M(2)).$$
By (\ref{A1cof}), since localization preserves cofibrations and $(A_1)_{K/2}=*$, there is an equivalence
\begin{equation}\label{MTM}\Sigma^4M_{K/2}\simeq T_{K/2}\w M(2),\end{equation} and hence
\begin{equation}\label{twosum}\pi_i(Y\w M(2))\approx\pi_{i-3}(M_{K/2})\oplus\pi_{i-4}(M_{K/2}).\end{equation}
 By \cite[4.2]{Bo79}, 
$$\pi_{8i+d}(M_{K/2})=\begin{cases}0&d=4,5\\
\Z/2&d=-2,3\\
\Z/2\oplus\Z/2&d=-1,2\\
\Z/4\oplus\Z/2&d=0,1.\end{cases}$$
This is the sum of two ``lightning flashes," one beginning in $8d-2$ and the other in $8d-1$.
Substituting this information into (\ref{twosum}) yields exactly the orders which were shown in the previous
paragraph to be true if and only if all the exact sequences (\ref{4seq}) split. \end{pf}

\section{Determination of $KO^*(\Phi_1DI(4))$}\label{Bosec}
In this section, we prove Theorem \ref{KO}, which shows how
$\psi^k$ in $K^*(BDI(4))$ leads to the determination of $KO^*(\Phi_1DI(4))$. 
Our presentation here follows suggestions in several e-mails from Pete Bousfield.

The first result explains how $KO^*(BDI(4))$ follows from $K^*(BDI(4))$.
\begin{thm}\label{And} There are classes $g_8$, $g_{12}$, and $g_{24}$ in
$KO^0(BDI(4))$ such that $c(g_i)=\xi_i$, with $\xi_i$ as in \ref{psik}, and
$$KO^*(BDI(4))\approx KO^*[\![g_8,g_{12},g_{24}]\!].$$
The Adams operations $\psi^2$ and $\psi^3$ mod decomposables on the basis of
$g_i$'s is as in \ref{psik}.\end{thm}
\begin{pf} In \cite[2.1]{And}, it is proved that if there is a torsion-free subgroup
$F^*\subset KO^*(X)$ such that $F^*\otimes K^*(pt)\to K^*(X)$ is an isomorphism, then
so is $F^*\otimes KO^*(pt)\to KO^*(X)$.
The proof
is a Five Lemma argument using exact sequences in \cite[p.257]{Kar}. 
Although the result is stated for ordinary (not 2-completed) $KO^*(-)$,
the same argument applies in the 2-completed context.  If $F^*$ is a multiplicative
subgroup, then the result holds as rings. Our result then follows from \ref{psik},
since the generators there are real.
 A similar proof can be derived from \cite[2.3]{Bo}.
\end{pf}

Next we need a similar sort of result about $KO^*(DI(4))$. We could derive much of
what we need by an argument similar to that just used, using the result of \cite{JO} 
about $K^*(DI(4))$ as input. However, as we will need this in a specific form in 
order to use it to draw conclusions about $KO^*(\Phi_1 DI(4))$, we begin by introducing
much terminology from \cite{Bo}.

The study of united $K$-theory begins with two categories, which will then be endowed
with additional structure. We begin with a partial definition of each, and their
relationship. For complete details, the reader will need to refer to \cite{Bo} or
an earlier paper of Bousfield.
\begin{defin}\label{CRdef} $(\cite[2.1]{Bo})$ A {\it $CR$-module} $M=\{M_C,M_R\}$
consists of $\Z$-graded $2$-profinite abelian groups $M_C$ and $M_R$ with continuous additive
operations $M_C^*@>B>\approx>M_C^{*-2}$, $M_R^*@>B_R>\approx>M_R^{*-8}$, $M_C^*@>t>\approx>M_C^*$, 
$M_R^*@>\eta>>M_R^{*-1}$, $M_R^*@>c>>M_C^*$, $M_C^*@>r>>M_R^*$, satisfying 15 relations, which we
will mention as needed.\end{defin}
We omit the descriptor ``2-adic," which Bousfield properly uses, just as we omit writing the 2-adic
coefficients $\hat\Z_2$ which are present in all our $K$- and $KO$-groups. 
\begin{expl} \label{CRexpl}For a spectrum or space $X$, the united $2$-adic $K$-cohomology
$$K_{CR}^*(X):=\{K^*(X),KO^*(X)\}$$
is a $CR$-module, with complex and real Bott periodicity, conjugation, the Hopf map, complexification,
and realification giving the respective operations.\end{expl}
\begin{defin}\label{deldef} $(\cite[4.1]{Bo})$ A $\Delta$-module $N=\{N_C,N_R,N_H\}$ is a triple of 
$2$-profinite abelian groups $N_C$, $N_R$, and $N_H$ with continuous additive operations
$N_C@>t>\approx>N_C$, $N_R@>c>>N_C$, $N_C@>r>>N_R$, $N_H@>c'>>N_C$, and $N_C@>q>>N_H$ satisfying
nine relations.\end{defin}
\begin{expl}\label{deltaexpl} For a $CR$-module $M$ and an integer $n$, there is a $\Delta$-module $\Delta^nM=\{M_C^n,M_R^n,M_R^{n-4}\}$
with $c'=B^{-2}c$ and $q=rB^2$. In particular, for a space $X$ and integer $n$, there is a $\Delta$-module
$K^n_{\Delta}(X):=\Delta^nK^*_{CR}(X)$.\end{expl}

Now we add additional structure to these definitions.
\begin{defin}\label{AdDeldef} $(\cite[4.3,6.1]{Bo})$ A $\theta\Delta$-module is a $\Delta$-module
$N$ together with homomorphisms $N_C@>\theta>>N_C$, $N_R@>\theta>>N_R$, and $N_H@>\theta>>N_R$
satisfying certain relations listed in \cite[4.3]{Bo}.
An {\it Adams $\Delta$-module} is a $\theta\Delta$-module $N$ together with
Adams operations $N@>\psi^k>\approx >N$ for odd $k$ satisfying the familiar properties.
\end{defin}
\begin{expl} In the notation of Example \ref{deltaexpl}, $K_{\Delta}^{-1}(X)$ is an Adams $\Delta$-module
with $\theta=-\lambda^2$.\end{expl}
\begin{defin} $(\cite[2.6,3.1,3.2]{Bo})$ A special $\phi CR$-algebra $\{A_C,A_R\}$ is a $CR$-module with
bilinear $A^m_C\times A_C^n\to A_C^{m+n}$ and $A^m_R\times A_R^n\to A_R^{m+n}$ and also $A_C^0@>\phi>> A_R^0$
and $A_C^{-1}@>\phi>> A_R^0$ satisfying numerous properties.\end{defin}
\begin{rmk} {\rm The operations $\phi$, which were initially defined in \cite{BoK}, are less familiar
than the others. Two properties are $c\phi a=t(a)a$ and $\phi(a+b)=\phi a+\phi b+r(t(a)b)$ for $a,b\in A_C^0$.
For a connected space $X$, $K_{CR}^*(X)$ is a special $\phi CR$-algebra.}\end{rmk}

The following important lemma is taken from \cite{Bo}.
\begin{lem} $(\cite[4.5,4.6]{Bo})$ For any $\theta\Delta$-module $M$, there is a universal
special $\phi CR$-algebra $\hat LM$. This means that there is a morphism $M@>\a>> \hat LM$ such that
any morphism from $M$ into a $\phi CR$-algebra  factors as $\a$ followed by a unique $\phi CR$-algebra morphism.
There is an algebra isomorphism $\hat\Lambda M_C\to(\hat LM)_C$, where $\hat\Lambda(-)$ is the 2-adic
exterior algebra functor.\end{lem}

In \cite[2.7]{Bo}, Bousfield defines, for a $CR$-algebra $A$, the indecomposable quotient
$\hat QA$. We apply this to $A=K^*_{CR}(BDI(4))$, and consider the $\Delta$-module $\hat QK^0_\Delta(BDI(4))$,
analogous to \cite[4.10]{Bo}. We need the following result, which is more delicate than the
$K^{-1}_\Delta$-case considered in \cite[4.10]{Bo}.
\begin{lem} With $\theta=-\lambda^2$, the $\Delta$-module $\hat QK^0_\Delta(BDI(4))$ becomes
a $\theta\Delta$-module.\end{lem}
\begin{pf} First we need that $\theta$ is an additive operation. In \cite[3.6]{BoK}, it is shown
that $\theta(x+y)=\theta(x)+\theta(y)-xy$ if $x,y\in KO^n(X)$ with $n\equiv0\mod4$. The additivity
follows since we are modding out the product terms. (In the case $n\equiv-1\mod 4$ considered in 
\cite[4.10]{Bo}, the additivity of $\theta$ is already present before modding out indecomposables.)

There are five additional properties which must be satisfied by $\theta$.
That $\theta cx=c\theta x$ and $\theta tz=t\theta z$ are easily obtained from \cite[3.4]{BoK}.
That $\theta c'y=c\theta y$ follows from \cite[6.2(iii),6.4]{BoTop}. That $\theta qz=\theta rz$ follows
from preceding \cite[3.10]{BoK} by $c$, which is surjective for us. Here we use that $rc=2$ and $q=rB^2$.
Finally, $\theta rz=r\theta z$ for us, since $c$ is surjective; here we have used the result $\bar\phi cx=0$
given in \cite[4.3]{Bo}.\end{pf}

Now we obtain the following important description of the $CR$-algebra $K^*_{CR}(DI(4))$.
\begin{thm}\label{DI4} There is a morphism of $\theta\Delta$-modules
$$\hat QK_\Delta^0(BDI(4))\to \tilde K_\Delta^{-1}(DI(4))$$
which induces an isomorphism of special $\phi CR$-algebras
$$\hat L(\hat QK_\Delta^0(BDI(4)))\to K^*_{CR}(DI(4)).$$
\end{thm}
\begin{pf} The map $\Sigma DI(4)=\Sigma\Omega BDI(4)\to BDI(4)$ induces a morphism
$$K_\Delta^0(BDI(4))\to K^{-1}_\Delta(DI(4))$$ which factors through the indecomposable quotient
$\hat QK_\Delta^0(BDI(4))$.
In \cite[1.2]{JO}, a general result is proved which implies that $K^*(DI(4))$ is an exterior algebra
on elements of $K^1(DI(4))$ which correspond to the generators of the power series algebra
$K^*(BDI(4))$ under the above morphism followed by the Bott map. Thus our result will follow
from \cite[4.9]{Bo}, once we have shown that the $\theta\Delta$-module $M:=\hat QK_\Delta^0(BDI(4))$
is robust(\cite[4.7]{Bo}).  This requires that $M$ is profinite, which follows as in the remark following
\cite[4.7]{Bo}, together with two properties regarding $\bar\phi$, where $\bar\phi z:=\theta rz-r\theta z$
for $z\in M_C$. In our case, $c$ is surjective, and so $\bar\phi=0$ as used in the previous proof.

One property is that $M$ is torsion-free and exact. This follows from the Bott exactness of
the $CR$-module $K_{CR}^*(BDI(4))$ noted in \cite[2.2]{Bo}, and \cite[5.4]{Bo}, which states that,
for any $n$, the $\Delta$-module $\Delta^nN$ associated to a Bott exact $CR$-module $N$ with $N^n_C$
torsion-free and $N^{n-1}_C=0$ is torsion-free. The other property is $\ker(\bar\phi)=cM_R+c'M_H+2M_C$.
For us, both sides equal $M_C$ since $c$ is surjective and $\bar\phi=0$.\end{pf}

Our Theorem \ref{KO} now follows from \cite[9.5]{Bo} once we have shown that the Adams $\Delta$-module
$M:=\hat QK_\Delta^0(BDI(4))$ is ``strong." (\cite[7.11]{Bo}) This result (\cite[9.5]{Bo}) requires that
the space (here $DI(4)$) be an $H$-space (actually $K/2_*$-durable, which is satisfied by $H$-spaces)
and that it satisfies the conclusion of our \ref{DI4}. It then deduces that $KO^*(\Phi_1DI(4))$ fits into an
exact sequence which reduces to ours provided $M_R=M_C$ and $M_H=2M_C$. These equalities are implied by 
$\hat QK_\Delta^0(BDI(4))$ being exact, as was noted to be true in the previous proof, plus $t=1$ and $c$ surjective, as were noted to be true in \ref{psik}. Indeed, the exactness
property, (\cite[4.2]{Bo}), includes that $cM_R+c'M_H=\ker(1-t)$ and $cM_R\cap c'M_H=\im(1+t)$. Another
perceptible difference is that Bousfield's exact sequence is in terms of $\bar M:=M/\bar\phi$, while ours involves
$M$, but these are equal since, as already observed, $\bar\phi=0$ since $c$ is surjective.

Note also that the Adams operations in $\bar M$ in the exact sequence of \cite[9.5]{Bo}, which reduces to that
in our \ref{KO}, are those in the Adams $\Delta$-module $\hat QK_\Delta^0(BDI(4))$, which are given
in our \ref{psik}. The morphism $\theta$ in \cite[9.5]{Bo} or our \ref{KO} is $\frac12\psi^2$, since
this equals $-\lambda^2$ mod decomposables.

Finally, we show that our $M$ is strong. One of the three criteria for being strong is to be robust,
and we have already discussed and verified this. The second requirement for an Adams $\Delta$-module
to be strong is that it be ``regular."  This rather technical condition is defined in \cite[7.8]{Bo}.
In \cite[7.9]{Bo}, a result is proved which immediately implies that $\tilde K_\Delta^{-1}(DI(4))$ is
regular. By  \ref{DI4}, our $M$ injects into
$\tilde K_\Delta^{-1}(DI(4))$, and so by \cite[7.10]{Bo}, which states that a submodule of a regular
module is regular, our $M$ is regular.

The third requirement for $M$ to be strong is that it be $\psi^3$-splittable (\cite[7.2]{Bo}), which
means that the quotient map $M\to M/\bar\phi$ has a right inverse. As we have noted several times,
we have $\bar\phi=0$, and so the identity map serves as a right inverse to the identity map.
This completes the proof that our $M$ is strong, and hence that \cite[9.5]{Bo} applies to $DI(4)$
to yield our Theorem \ref{KO}.

\def\line{\rule{.6in}{.6pt}}

\end{document}